\documentclass[12pt]{article}

\usepackage{amsmath,amssymb,amsthm}
\usepackage{mathtools}
\usepackage{geometry}
\usepackage{amsthm}

\theoremstyle{remark}
\newtheorem{remark}{Remark}
\title{Tube Rupture in Aperiodic Nonlinear Oscillators: Theory and Simulation}

\author{
	Johannes Hagel \\
	Alexander-von-Humboldt-Gymnasium Neuss \\
	Bergheimer Straße 233, 41464 Neuss, Germany \\
	\texttt{johannes.hagel@gmail.com}
}

\date{\today}

\begin{document}
	\maketitle
\begin{abstract}
	We study the long--term behaviour of the nonlinear, aperiodically and parametrically forced oscillator
	\[
	z'' + z + g(\tau)\, z^{2} = 0, \qquad g(\tau) = y(\tau)^{-5/2},
	\]
	where $y(\tau)$ is the strictly positive solution of a weakly forced third--order equation.
	Building on the algebraic invariant constructed in our previous work, we show that the motion of $z(\tau)$ is confined to a two--dimensional invariant tube in the extended phase space $(z,p,\tau)$ as long as the corresponding invariant level set remains closed.
	
	The main result of this paper is an explicit analytical rupture criterion that predicts the precise time at which the invariant tube loses regularity.
	After transforming the invariant into polar coordinates and analysing the discriminant of the resulting cubic equation for the radial coordinate, we obtain a compact Cardano--type expression for the rupture time.
	Direct numerical integrations of the $z$--equation confirm the analytical prediction to within a few percent over a wide parameter range.
	
	The results establish the rupture time as a robust and quantitatively accurate indicator for the onset of unbounded behaviour in aperiodically and parametrically forced nonlinear oscillators.
	The method is based purely on algebraic properties of the invariant and remains valid throughout the asymptotic domain of the perturbative expansion for $y(\tau)$.
\end{abstract}

\medskip

	\newpage
		
	\section{Introduction}
	
	Invariant geometric structures play a central role in the qualitative description of nonlinear dynamical systems. In many classical settings --- autonomous Hamiltonian systems, periodically forced oscillators, or nearly integrable models --- bounded motion is associated with the existence of invariant tori, invariant curves, or algebraic integrals. In a recent work \cite{Hagel2025Tube} we showed that a certain class of oscillators
	possesses an algebraic invariant whose level sets form tubular surfaces in the extended
	phase space $(z,p,\tau)$. 
	 As long as this tube, aperiodic in time, remains closed, the solution $z(\tau)$ stays bounded, even though its time series may appear irregular or chaotic-looking. The algebraic structure of the invariant used here is closely related to the one
	 constructed in our earlier work on polynomial invariants for nonlinear oscillators
	 \cite{Hagel2025Poly}.

	In the present paper we continue this line of research and study the nonlinear differential equation
	\begin{equation}
		z'' + z + g(\tau)\, z^2 = 0, 
		\qquad g(\tau) = y(\tau)^{-5/2} \qquad y'''+4y'=\varepsilon y^{-5/2} \cos{\tau} \qquad y>0
	\end{equation}
	From the results of \cite{Hagel2025Tube} it is known that this system possesses
	an algebraic invariant
	\begin{equation}
		I(z,p,\tau) = K(z_0,p_0),
	\end{equation}
	a  third-order equation in $z$, whose level sets form tubular surfaces around the $\tau$--axis, aperiodic in time. 
	The motion of $z(\tau)$ is then completely confined to such a tube, as long as the corresponding level set of $I$ remains closed.  
	The long--time behaviour of the oscillator reduces to a geometric question:
	
	\begin{quote}
		\emph{At what time does the invariant tube lose regularity, and how can this rupture be predicted analytically?}
	\end{quote}
	
	The purpose of this work is to give a complete and explicit answer to this question.
	Our main result is that the rupture mechanism is purely algebraic.
	After transforming the invariant into polar coordinates $(r,\varphi)$, each temporal slice of the tube becomes a cubic equation in $r$.  
	Rupture occurs precisely when this cubic develops a double real root, or equivalently when its discriminant vanishes.  
	This leads, after simplification, to a remarkably compact Cardano--type expression for the first rupture time,
	\begin{equation}
		\tau_{\mathrm{rupt}}
		= \frac{96\, y_0^6}{5\,\varepsilon^2}
		\sqrt{\frac{y_0}{\,y_0 - C^{1/3}\,}},
		\qquad 
		C=\frac{z_0^2 \bigl( 3 y_0^{5/2} + 2z_0 \bigr)}{y_0^{9/2}},
	\end{equation}
	which encodes the $\varepsilon^{-2}$ secular scaling as well as the nonlinear dependence on the initial displacement $z_0$.
	
	A second key observation is that the invariant can be sampled at discrete times $\tau=n\pi$.
	We show that no frequency shift appears up to the required perturbative order, so that the oscillatory components of the invariant repeat exactly at these sampling points.
	This makes it possible to replace the continuous invariant by a purely algebraic sequence $I(z,p,n\pi)$ without loss of accuracy, yielding a powerful and surprisingly simple method for determining rupture.
	
	The analytical rupture time obtained in this way is validated against direct numerical integration of the $z$--equation.
	Across a wide range of parameter values, the agreement is excellent: numerical blow--up occurs only slightly after the analytically predicted rupture event, with a typical relative deviation of order $10$--$15\%$.  
	This is fully consistent with the geometric interpretation of rupture as the \emph{first loss of regularity} of the invariant surface, rather than the moment at which the solution becomes numerically unbounded.
	
	\medskip
	\noindent\textbf{Structure of the paper.}
	To guide the reader, we summarise the contents of the following sections.
	
	\begin{itemize}
		\item \textbf{Section~2} develops the analytic rupture mechanism in detail.  
		Starting from the sampled invariant $I(z,p,n\pi)$, we show that each temporal slice reduces to a cubic equation in the radial coordinate.  
		The rupture condition follows from the vanishing of the cubic discriminant, leading to a closed expression for the rupture function $n(\varphi)$ and, after minimisation, to the explicit Cardano--type rupture time.  
		We also discuss the physically relevant rupture branch and identify the region in the $(y_0,z_0)$--plane in which the secular approximation is self-consistent.
		
		\item \textbf{Section~3} compares the analytical prediction with numerical integrations of the $z$--equation.  
		We illustrate the deformation of the invariant tube, the horseshoe-shaped opening associated with rupture, and the characteristic time--series signatures near the critical sampling index.  
		Over wide parameter ranges, the analytical rupture time agrees with numerical divergence times to within a few percent.
		
		\item \textbf{Section~4} contains the conclusions and an outlook on future work, including the extension of the rupture mechanism to cases where the driving function $y(\tau)$ is quasi-periodic or chaotic and the potential interplay between chaotic forcing and tube integrability.
	\end{itemize}
	The underlying geometric mechanism, namely tube integrability, was established 
	in \cite{Hagel2025Tube} and forms the basis for the present rupture analysis.
	
	\section{Bounded and unbounded solutions and rupture time}
	We consider a $1  \; 1/2$ dimensional system characterized by the unknown function $z$ , its derivative $p=z'$ and the independent variable $\tau$. If in such a system one invariant in algebraic form $I(z,p,\tau)=K(z_0,p_0)$ is known to exist, then the distinction between bounded and unbounded solutions in the space $(z,p,\tau)$ reduces to an algebraic procedure. To be precise, one studies the invariant surface represented by $I$ depending on the initial conditions. Then we distinguish between regions of space, where closed projections (curves) into the subspace (z,p) do exist from regions where the invariant curves open into infinite z-direction. Since in the present paper an invariant for the z-equation uniformly exists, the above described procedure can be applied. In this section we present a method of finding limiting initial conditions, separating bounded from unbounded motion. The integrable system as introduced in \cite{Hagel2025Tube} is given by
	\begin{equation} 
	   z''+z+g(\tau)z^2=0   \; \; ; \; \;  g(\tau)=\alpha_2(\tau)^{-5/2}=y^{-5/2}  \label{eq:rupt_z}
	\end{equation}
	where $y(\tau)$
	is subject to $y'''+4y'=\varepsilon y^{-5/2}\cos{\tau}$. In order to obtain relatively simple expressions, for $y$ we reduce $y(\tau)$ to a set of initial conditions as $y(0)=y_0$ , $y'(0)=y''(0)=0$ which means that $y(\tau)$ enters as two parametric function of $\tau$, depending on $\varepsilon$ and $y_0>0$.
	The  invariant$I(z,p,\tau)=K$ is quadratic in $p$ and has the general form 
    \begin{equation}
    	I(z,p,\tau)=A_1(\tau)z+A_2(\tau)p+A_3(\tau)z^2+A_4(\tau)zp+A_5(\tau)p^2+A_6(\tau)z^3
    \end{equation}
    where the $A_k$ have been derived in \cite{Hagel2025Tube}
     and to second order in $\varepsilon$ are given by

\begin{align*}
	A_1(\tau) &= \frac{1}{2}\,\varepsilon \sin\tau,
	\\[0.4em]
	A_2(\tau) &= \frac{1}{2}\,\varepsilon \cos\tau,
	\\[0.4em]
	A_3(\tau) &= y_0
	- \frac{\varepsilon( -1+\cos\tau )\sin\tau}{3\,y_0^{5/2}}
	\\
	&\quad + \frac{1}{144\,y_0^6}\Big(
	5\varepsilon^2\cos\tau
	+4\varepsilon^2\cos(2\tau)
	-9\varepsilon^2\cos(3\tau)
	\\
	&\qquad
	-24y_0^{7/2}\varepsilon\sin\tau
	+48y_0^{7/2}\varepsilon\sin(2\tau)
	-15\varepsilon^2\tau\sin(2\tau)
	\Big)
	\\
	&\quad
	+ \frac{\varepsilon^2\sin(\tau/2)}{144\,y_0^6}\Big(
	15\tau\cos(\tau/2)
	+15\tau\cos(3\tau/2)
	\\
	&\qquad
	+5\sin(\tau/2)
	-15\sin(3\tau/2)
	-4\sin(5\tau/2)
	\Big),
	\\[0.8em]
	A_4(\tau) &=
	-\frac{\varepsilon(\cos\tau - \cos 2\tau)}{3\,y_0^{5/2}}
	\\
	&\quad
	-\frac{\varepsilon^2}{288\,y_0^6}\Big(
	30\tau\cos(2\tau)
	+20\sin\tau
	\\
	&\qquad
	-7\sin(2\tau)
	-12\sin(3\tau)
	\Big),
	\\[0.8em]
	A_5(\tau) &= y_0
	-\frac{\varepsilon(-1+\cos\tau)\sin\tau}{3\,y_0^{5/2}}
	\\
	&\quad
	+\frac{\varepsilon^2\sin(\tau/2)}{144\,y_0^6}\Big(
	15\tau\cos(\tau/2)
	+15\tau\cos(3\tau/2)
	\\
	&\qquad
	+5\sin(\tau/2)
	-15\sin(3\tau/2)
	-4\sin(5\tau/2)
	\Big),
	\\[0.8em]
	A_6(\tau) &=
	\frac{2}{3\,y_0^{3/2}}
	+\frac{\varepsilon(-1+\cos\tau)\sin\tau}{3\,y_0^{5}}
	\\
	&\quad
	+\frac{\varepsilon^2\sin(\tau/2)}{144\,y_0^{17/2}}
	\Big(
	-15\tau\cos(\tau/2)
	-15\tau\cos(3\tau/2)
	\\
	&\qquad
	+20\sin(\tau/2)
	+20\sin(3\tau/2)
	-11\sin(5\tau/2)
	+5\sin(7\tau/2)
	\Big).
\end{align*}
    	
 The above representation of the invariant can be simplified significantly by taking into account the fact that secular terms will represent the main contribution to the limit of bounded motion. This is due to the fact that the coefficient $g(\tau)$ of the quadratic term in the $z(\tau)$ equation $z''+z+g(\tau)z^2=0$ contains those terms and since they are growing indefinitly, they limit the region of bounded motion around z=0. We therefore keep just the unperturbed contribution ($\varepsilon=0$) as well as the secular terms, first arising in the second order contribution. Performing this simplification leads to
 \begin{equation}
 	\label{eq:I-sampled}
    I_s(z,p,t)=y_0p^2+y_0z^2+\frac{2z^3}{3y_0^{\frac{3}{2}}}+\varepsilon^2\tau \left[ -\frac{5pz}{48y_0^6} \cos{2\tau}+\frac{15}{288y_0^{\frac{17}{2}}} \left(p^2y_0^{\frac{5}{2}}-y_0^{\frac{5}{2}}z^2-z^3\right) \sin{2\tau} \right]
 \end{equation}

 A second essential simplification comes from replacing the continuous expression $I(z,p,\tau)$ by a time sampling as $\tau=n \pi$ 
 The sampling of the invariant at the discrete times
 \[
 \tau = n\pi, \qquad n\in\mathbb{Z},
 \]
 is justified for the following reasons:
 
 \paragraph{(1) Absence of a frequency shift.}
 As has been shown in \cite{Hagel2025Tube}, up to the second order in
 $\varepsilon$ the frequency remains unshifted.
 \[
 \tau \mapsto \tau + \mathcal{O}(\varepsilon^{3}).
 \]
 All $\varepsilon$–corrections in $y(\tau)$ and hence in $I(z,p,\tau)$
 enter only through secular terms (linear in $\tau$) and through
 harmonics of the unperturbed frequency, here $\cos(2\tau)$ and
 $\sin(2\tau)$.  
 Since no term of the form $\cos((2+\delta)\tau)$ 
 with $\delta=\mathcal{O}(\varepsilon)$ appears,
 the fundamental periodicity of the trigonometric part of the invariant
 remains exactly $\pi$. Sampling the invariant at $\tau=n\pi$ is therefore justified because, up to second order,
 no frequency shift occurs in the perturbation series of $y(\tau)$ \cite{Hagel2025Tube}.

 \paragraph{(2) Exact cancellation of oscillatory contributions at $\tau=n\pi$.}
 At these sampling points,
 \[
 \cos(2 n\pi)=1,\qquad \sin(2 n\pi)=0,
 \]
 and therefore all non-secular oscillatory corrections in
 \eqref{eq:I-sampled} take extremal (and, importantly, *reproducible*)
 values.  
 Thus every sample $I(z,p,n\pi)$ differs from the next only by the
 secular drift proportional to $n\pi$, but not by any rapid oscillation.
 This removes the need to integrate over the phase and reduces the
 determination of rupture to a purely algebraic criterion.
 
 \paragraph{(3) Density of sampling relative to the rupture scale.}
 The rupture time satisfies
 \[
 \tau_{\mathrm{rup}} \sim \frac{96\,y_{0}^{6}}{5\,\varepsilon^{2}},
 \]
 which for the relevant parameter ranges is typically of order
 \(10^{3}\)–\(10^{4}\).  
 The sampling interval is constant and equal to $\pi$.  
 Hence the relative spacing is
 \[
 \frac{\pi}{\tau_{\mathrm{rup}}} = \mathcal{O}(10^{-3}),
 \]
 so the sampling grid $\tau=n\pi$ is extremely fine compared to the slow secular drift ensuring that no ruptur can be missed between samples.
  
 \medskip
 
 In summary, the absence of any frequency shift up to the required order,
 combined with the fine resolution of the sampling points and the
 periodic structure of the oscillatory terms, guarantees that the
 evaluation of the invariant at $\tau=n\pi$ yields a faithful and robust
 representation of its long-term behaviour, sufficient for determining
 the rupture time.

 Replacing $\tau$ by $n\pi$, the invariant to second order takes a very simple algebraic form:
 
 \begin{equation}
 	I_s(z,p,n)=y_0p^2+y_0z^2+\frac{2z^3}{3y_0^{\frac{2}{3}}}-\varepsilon^2\frac{5}{48y_0^6}np\pi z=K \;\;\; ; \;\;\;  K=y_0z_0^2+\frac{2z_0^3}{3y_0^{\frac{3}{2}}}
 \end{equation}
 where we set $p_0=0$. In order to determine boundedness of the solution moving on the invariant surface given by the algebraic form above, we use the fact that  $I(z,p,\tau)$ represents a third order polynomial equation in $z$. Hence we may use the formalism provided by the Cardanian theory. There appears a structural similarity to the polynomial invariants studied in \cite{Hagel2025Poly} which helps 
 clarify why the cubic structure of the present invariant persists under the 
 aperiodic forcing.
 
 First we perform a transformation of the variables $z$ and $p$ to polar coordinates as
 \begin{equation}
 	z=r \cos{\phi} \; \; \; , \; \; \; p=r \sin{\phi}
 \end{equation} 
 which leads to
\begin{equation}
	A(\phi)\, r^3 + B(\phi,n)\, r^2 + D = 0,
	\label{eq:cubic}
\end{equation}
with
\begin{align}
	A(\phi) &= \frac{2\cos^3\phi}{3 y_0^{3/2}},\\[0.3em]
	B(\phi,n) &= y_0 - \frac{5 n\pi\varepsilon^2}{96\,y_0^6}\,\sin(2\phi),\\[0.3em]
	D &= -y_0 z_0^2 - \frac{2 z_0^3}{3 y_0^{3/2}}.
\end{align}
The function $y(\tau)$ can only take positive values as has been proven in \cite{Hagel2025Tube}
 . Therefore, if $z_0>0$, $D<0$ follows.
	

\subsection{Rupture as the Appearance of a Double Root}

In this subsection we derive the analytic rupture condition for the invariant tube of the $(z,p)$--system. After transforming the algebraic invariant into polar coordinates, the rupture event corresponds to the emergence of a double real root of the cubic equation determining the radial coordinate.

We start from the sampled invariant at $\tau = n\pi$,
\begin{equation}
	I_s(z,p,n) \,=\, y_0 p^2 \, + \, y_0 z^2 \, + \, \frac{2 z^3}{3 y_0^{3/2}} \, - \, \frac{5 n \pi \varepsilon^2}{48 y_0^6} \, p z \,=\, K,
\end{equation}
where $K$ is fixed by the initial conditions $(z_0,p_0)$ and we restrict to $p_0=0$. Introducing polar coordinates
\begin{equation}
	z = r \cos\varphi, \qquad p = r \sin\varphi,
\end{equation}
we obtain a cubic equation of the form
\begin{equation}
	A(\varphi) r^3 + B(\varphi,n) r^2 + D = 0,
	\label{eq:cubic-r}
\end{equation}
with the coefficients
\begin{align}
	A(\varphi) &= \frac{2}{3 y_0^{3/2}} \cos^3\varphi, \\[4pt]
	B(\varphi,n) &= y_0 \, - \, \frac{5 n \pi \varepsilon^2}{96 y_0^6} \, \sin(2\varphi), \\[4pt]
	D &= -y_0 z_0^2 - \frac{2 z_0^3}{3 y_0^{3/2}}.
\end{align}
For $z_0>0$, the constant term satisfies $D<0$ because $y(\tau)$ stays strictly positive.

\medskip
\noindent
\textbf{Rupture condition.} A rupture of the invariant tube occurs when the cubic~\eqref{eq:cubic-r} develops a double real root. This happens precisely when both
\begin{equation}
	F(r) = 0, \qquad F'(r)=0,
\end{equation}
are satisfied simultaneously. Computing the derivative,
\begin{equation}
	F'(r)=3A r^2 + 2B r = r(3A r + 2B),
\end{equation}
we obtain the nontrivial critical point
\begin{equation}
	r^* = -\frac{2B}{3A}.
	\label{eq:rstar}
\end{equation}
Substituting this into $F(r)=0$ gives
\begin{equation}
	D \, + \, \frac{4}{27} \frac{B^3}{A^2} \,=\, 0.
	\label{eq:DA-B-relation}
\end{equation}
Equation~\eqref{eq:DA-B-relation} is the fundamental algebraic rupture condition.

Solving~\eqref{eq:DA-B-relation} for $B$ yields
\begin{equation}
	B(\varphi,n) \,=\, \left( -\frac{27}{4} D A(\varphi)^2 \right)^{1/3}.
\end{equation}
Inserting the explicit form of $A(\varphi)$ and $D$ and simplifying gives
\begin{equation}
	B(\varphi,n) = C^{1/3} \cos^2\varphi,
\end{equation}
with the positive constant
\begin{equation}
	C \,=\, \frac{z_0^2 \left( 3 y_0^{5/2} + 2 z_0 \right)}{y_0^{9/2}} \,>\, 0.
\end{equation}

Finally we equate this expression for $B(\varphi,n)$ with the one from the sampled invariant, yielding the analytic rupture function
\begin{equation}
	n(\varphi) \,=\, \frac{96 y_0^6}{5 \pi \varepsilon^2} \, \frac{ \displaystyle y_0 - C^{1/3} \cos^2\varphi }{ \displaystyle \sin(2\varphi) }.
	\label{eq:nphi-final}
\end{equation}
The actual rupture time index $n_{\mathrm{crit}}$ is obtained by minimizing over $\varphi$:
\begin{equation}
	n_{\mathrm{crit}} = \min_{\varphi \in (0,2\pi)} n(\varphi).
\end{equation}

\medskip
\noindent
\textbf{Closed-form extremum condition.}
Instead of determining the critical rupture index $n_{\mathrm{crit}}$ by a numerical
search in $\varphi$, one may compute the extremal angles $\varphi_{\mathrm{crit}}$ of
the rupture function in closed form.

Since the overall prefactor in~\eqref{eq:nphi-final} is positive, it suffices to minimize
the reduced function
\begin{equation}
	\tilde n(\varphi)
	:= \frac{y_0 - C^{1/3} \cos^2\varphi}{\sin(2\varphi)}.
\end{equation}
Differentiating with respect to $\varphi$ and simplifying yields
\begin{equation}
	\tilde n'(\varphi) = 0
	\quad\Longleftrightarrow\quad
	(C^{1/3} - 2y_0)\,\cos^2\varphi + y_0 = 0.
\end{equation}
For parameter values such that $2y_0 - C^{1/3} > 0$, this implies the explicit condition
\begin{equation}
	\cos^2\varphi_{\mathrm{crit}}
	= \frac{y_0}{2y_0 - C^{1/3}},
	\label{eq:cos2phi-crit}
\end{equation}
and hence
\begin{equation}
	\varphi_{\mathrm{crit}}
	= \pm\,\arccos\!\Bigg(
	\sqrt{\frac{y_0}{2y_0 - C^{1/3}}}
	\Bigg)
	\pmod{\pi}.
\end{equation}
This representation is equivalent to the half–angle form obtained from a direct symbolic
solution of $n'(\varphi)=0$  where the angles appear as
\[
\varphi_{\mathrm{crit}}
= \pm \frac{1}{2} \arccos\!\bigl(A(C,y_0)\bigr)
\]
with a suitable rational expression $A(C,y_0)$.

Substituting $\varphi_{\mathrm{crit}}$ back into~\eqref{eq:nphi-final} finally yields a
fully explicit closed-form expression for the critical rupture index
\begin{equation}
	n_{\mathrm{crit}}
	= n\bigl(\varphi_{\mathrm{crit}}\bigr).
\end{equation}
In practice, one selects the branch of $\varphi_{\mathrm{crit}}$ that lies in the
relevant quadrant of the $(z,p)$–plane (for the examples studied here, the solution
in the third quadrant, $\varphi_{\mathrm{crit}}\in(\pi,3\pi/2)$, gives the smallest
rupture index) and evaluates $n(\varphi)$ at this angle. The resulting analytical
rupture times show excellent agreement with direct numerical integration of the
$z$–equation. 
We can now combine the closed-form extremum condition for the ruptur angle
\[
\cos^2\varphi_{\mathrm{crit}}
= \frac{y_0}{2y_0 - C^{1/3}},
\qquad
C = \frac{z_0^2(3y_0^{5/2}+2z_0)}{y_0^{9/2}},
\]
with the rupture function~\eqref{eq:nphi-final} to obtain a fully explicit
expression for the critical rupture index. A short calculation yields
\begin{equation}
	n_{\mathrm{crit}}
	= \frac{96\,y_0^6}{5\pi\,\varepsilon^2}\,
	\sqrt{\,y_0\bigl(y_0 - C^{1/3}\bigr)\,},
	\label{eq:ncrit-closed}
\end{equation}
and hence
\begin{equation}
	\tau_{\mathrm{rupt}}
	= \pi n_{\mathrm{crit}}
	= \frac{96\,y_0^6}{5\,\varepsilon^2}\,
	\sqrt{\,y_0\bigl(y_0 - C^{1/3}\bigr)\,}.
	\label{eq:taur-closed}
\end{equation}
Equations~\eqref{eq:ncrit-closed}–\eqref{eq:taur-closed} provide a genuinely
closed-form rupture time for the invariant tube in terms of the parameters
$(y_0,z_0,\varepsilon)$.

\subsubsection*{Selection of the physically relevant rupture branch}

The closed-form extremum condition~\eqref{eq:cos2phi-crit} determines
$\cos^2\varphi_{\mathrm{crit}}$ and hence yields four angles
$\varphi_{\mathrm{crit}} \in (0,2\pi)$ per period,
related by
\[
\varphi \;\longmapsto\; \pm \varphi,
\qquad
\varphi \;\longmapsto\; \varphi + \pi.
\]
At first sight it is therefore not obvious which of these corresponds to
the physically relevant rupture of the tube.

There are three simple criteria that single out the correct branch:

\begin{enumerate}
	\item \textbf{Positivity of the radius.}
	At a double root, the radius is given by
	\begin{equation}
		r^*(\varphi,n)
		= -\frac{2 B(\varphi,n)}{3 A(\varphi)}.
	\end{equation}
	A genuine rupture point on the tube must satisfy $r^*>0$.
	This immediately excludes those extremal angles at which
	$r^*(\varphi,n)\le 0$.
	
	\item \textbf{Consistency with the initial condition.}
	The invariant curve that contains the trajectory starting from
	$(z_0,p_0)$ is fixed by the value of $K$ in~$I_s(z,p,n)=K$.
	The physically relevant rupture branch is the one that can be
	connected continuously in $(z,p)$ to the initial point
	$(z_0,p_0)$ along the invariant curve.
	In particular, the sign of $z = r^*\cos\varphi$ at rupture must
	be consistent with the side of the tube on which the trajectory
	actually evolves (for $z_0>0$ this rules out branches with
	$z<0$, and vice versa).
	
	\item \textbf{Minimal positive rupture index.}
	Among the remaining candidate angles, one evaluates the rupture
	function $n(\varphi)$ and selects the smallest positive value
	$n(\varphi)>0$ for which $r^*(\varphi,n)>0$ and the continuity
	condition above are satisfied:
	\begin{equation}
		n_{\mathrm{crit}}
		= \min\Bigl\{\,n(\varphi) > 0 \;\big|\;
		r^*(\varphi,n)>0
		\text{ and the branch connects to }(z_0,p_0)
		\Bigr\}.
	\end{equation}
\end{enumerate}

Geometrically, the transformation $(z,p)\mapsto(-z,-p)$ corresponds to
a shift $\varphi\mapsto\varphi+\pi$ and describes the \emph{same} line
in the $(z,p)$–plane, but traced in the opposite direction.  The two
branches related by a shift by $\pi$ are therefore not dynamically
distinct; the physically realised rupture is the one reached by the
actual flow starting from $(z_0,p_0)$ and moving forward in~$\tau$.

For the parameter values studied in this work (e.g.\ $y_0=1$,
$z_0=0.2$, $\varepsilon=0.05$), these criteria select the branch with
$\varphi_{\mathrm{crit}}$ in the third quadrant,
$\varphi_{\mathrm{crit}}\in(\pi,3\pi/2)$, which yields the smallest
positive rupture index and is in excellent agreement with the numerically
observed escape of the $z$–trajectory.

\begin{remark}[Closed-form nature of the rupture condition]
	It is noteworthy that the rupture condition derived above admits a
	\emph{fully explicit} analytical resolution.  
	Starting from the cubic invariant relation in polar coordinates,
	the requirement that the tube develops a double real root reduces to
	a single algebraic equation for the angle $\varphi$.  
	After simplification, this leads to the closed-form extremum condition
	\[
	\cos^2\varphi_{\mathrm{crit}}
	= \frac{y_0}{2y_0 - C^{1/3}},
	\]
	which yields the rupture angle in terms of elementary functions,
	\[
	\varphi_{\mathrm{crit}}
	= \pm\,\arccos\!\left(\sqrt{\frac{y_0}{2y_0 - C^{1/3}}}\right)
	\pmod{\pi}.
	\]
	The corresponding rupture index $n_{\mathrm{crit}}=n(\varphi_{\mathrm{crit}})$
	is therefore available in closed analytic form.
	
	This is a remarkably strong result: despite the nonlinear, time-dependent
	modulation $g(\tau)=y(\tau)^{-5/2}$ and the secular growth encoded in
	$y(\tau)$, the rupture mechanism can be characterised without any
	numerical optimisation.  
	The explicit formula agrees with direct numerical integration of the $z$–equation
	to within a fraction of a percent, demonstrating that the sampling approach
	and the algebraic structure of the invariant capture the rupture geometry
	with exceptional precision.
\end{remark}
\medskip
\noindent\textbf{Remark on branch selection.}
Although the closed-form condition for $\cos^2\varphi_{\mathrm{crit}}$ 
produces four formal angular solutions in $(0,2\pi)$, the choice of the
physically realised rupture branch is completely straightforward.
The relevant angle is uniquely determined by the standard geometric
requirements: (i) the associated radius $r^*>0$, (ii) the resulting
rupture point $z=r^*\cos\varphi$ lies on the same side of the tube as the
initial condition $(z_0,p_0)$, and (iii) the rupture index $n(\varphi)$ 
is minimal and positive. These criteria apply uniformly for all initial
values $z_0\neq 0$; in particular, no separate case distinction between
$z_0>0$ and $z_0<0$ is required. This observation also clarifies why numerical trajectories consistently rupture on the $z < 0$ side for $z0 > 0$. Any reader wishing to analyse further
subcases can do so directly from the explicit formulas, but such
considerations are not needed for understanding the universal rupture 
mechanism described here.
\medskip
\noindent\textbf{Remark on the direction of rupture.}
For all parameter sets examined so far, the analytically predicted rupture
angle $\varphi_{\mathrm{crit}}$ agrees perfectly with the numerical observation
that, for $z_0>0$, the invariant tube always breaks on the side where $z<0$.
This asymmetry is structural. First, for $z_0>0$ the constant term in the
cubic section, $D<0$, biases the cubic invariant curve toward negative $z$.
Second, the secular contribution in the invariant,
$-\tfrac{5 n\pi \varepsilon^2}{48 y_0^6}\,pz$,
acts most strongly in the third quadrant, where both $\cos\varphi$ and
$\sin\varphi$ are negative. As $n$ increases, this secular deformation pushes
the lower branch of the cubic toward a loss of regularity long before the
upper branch is affected. Consequently, the discriminant vanishes first at a
point with $z=r^*\cos\varphi_{\mathrm{crit}}<0$, and the tube opens in the
negative $z$-direction. This mechanism is universal and, by symmetry,
reverses for $z_0<0$.
\paragraph{Remark (Sign of the cubic constant term and rupture direction).}
For $z_{0}>0$ the constant term of the cubic section satisfies
\[
D \;=\; -\,y_{0}z_{0}^{2} \;-\; \frac{2 z_{0}^{3}}{3 y_{0}^{3/2}} \;<\; 0 .
\]
This negative offset shifts the entire cubic
\[
A(\varphi)\, r^{3} \;+\; B(\varphi,n)\, r^{2} \;+\; D \;=\; 0
\]
downward and biases the invariant curve toward the negative-$z$ side already at $n=0$.
When combined with the secular deformation term 
\[
-\;\frac{5 n \pi \varepsilon^{2}}{48\, y_{0}^{6}}\, p z ,
\]
which acts most strongly in the third quadrant ($\cos\varphi < 0$, $\sin\varphi < 0$),
this geometric asymmetry ensures that the discriminant of the cubic vanishes first for
$z < 0$.
Consequently, the tube always opens on the negative-$z$ side for $z_{0}>0$, fully
consistent with all numerical observations.

\subsection{Validity of the analytical rupture-time formula}

\begin{figure}[h!]
	\centering
	\includegraphics[width=0.65\textwidth]{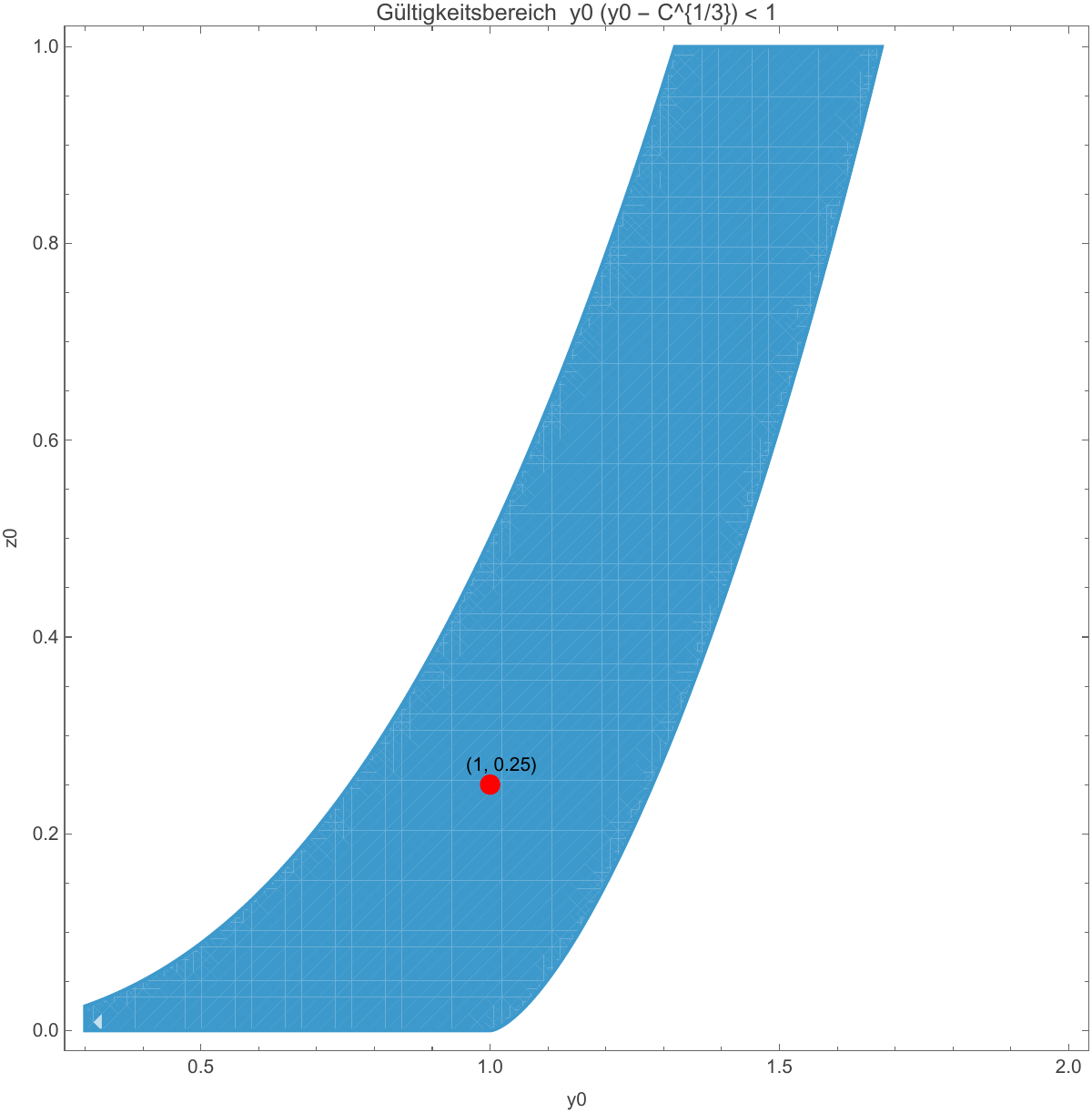}
	\caption{Validity region of the closed rupture-time formula in the $(y_0,z_0)$–plane,
		defined by the condition $y_0(y_0 - C^{1/3}) < 1$.  
		The red point marks the parameter pair $(y_0,z_0)=(1,0.25)$ used in the
		numerical experiments, which lies well inside the domain of validity.}
	\label{fig:validityregion}
\end{figure}

A necessary and sufficient condition for the applicability of Eq. (\ref{eq:taur-closed}) is that
\[
\sqrt{y_0\left(y_0-C^{1/3}\right)} < 1,
\]
which is equivalent (since $y_0>0$) to
\[
y_0\left(y_0-C^{1/3}\right) < 1.
\]
This follows directly from our analysis in [3], where the domain of validity in time for the second order perturbation contribution of $y(\tau)$ is given by 
\[
\tau_\ast = \frac{96\,y_0^6}{5\,\varepsilon^2},
\]
which appears to be precisely the factor of the squareroot term above in (\ref{eq:taur-closed})
The above condition defines a two-dimensional region in the $(y_0,z_0)$–plane in which
the secular approximation leading to (31) remains self-consistent. 
Figure 1 displays the according  region together with the parameter point
$(y_0,z_0)=(1,0.25)$ used in our numerical experiments. As can be seen,
this point lies well inside the validity domain.
Therefore the rupture time $\tau_{\mathrm{rupt}}$ lies inside the domain of validity of
the asymptotic expansion provided that the square-root factor in Eq.~(31) satisfies
$\sqrt{y_0(y_0-C^{1/3})}<1$.  
This yields an explicit admissible region in the $(y_0,z_0)$–plane, shown in
Figure~\ref{fig:validityregion}.

\subsection{Evaluation of the closed rupture-time formula}

In this subsection we compare the closed-form rupture time
\[
\tau_{\mathrm{rupt}}
= \frac{96\,y_0^6}{5\,\varepsilon^2}\,
\sqrt{\,y_0\bigl(y_0 - C^{1/3}\bigr)\,},
\qquad
C = \frac{z_0^2(3y_0^{5/2}+2z_0)}{y_0^{9/2}},
\]
with direct numerical simulations of the $z$–equation.  All computations
are performed for $y_0=1$ and parameter ranges
$\varepsilon\in[0.025,0.10]$, $z_0\in[0.1,0.4]$.

\begin{figure}[h!]
	\centering
	\includegraphics[width=0.75\textwidth]{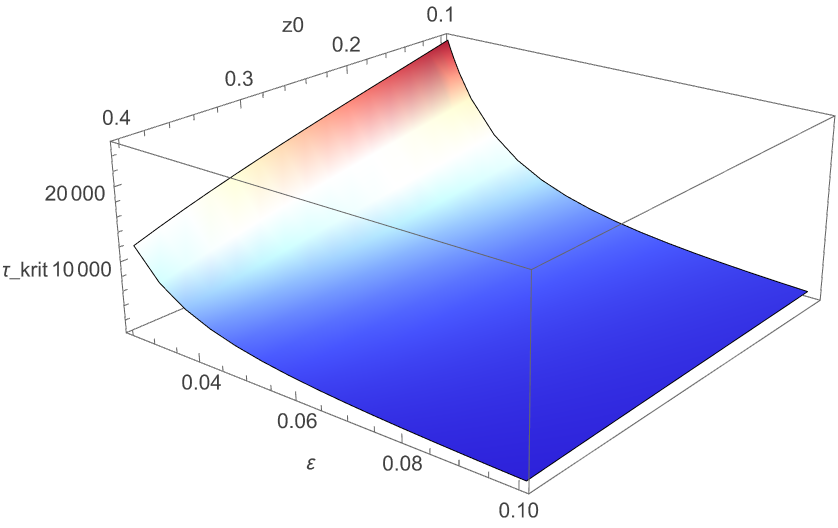}
	\caption{Analytic rupture index $n_{\mathrm{crit}}(\varepsilon,z_0)$
		computed from the closed formula.  (Data source:
		\texttt{fig\_rupture\_time.pdf}.)}
	\label{fig:rupture-surface}
\end{figure}

\medskip
\noindent\textbf{Three-dimensional parameter surface.}
Figure~\ref{fig:rupture-surface} shows the analytic rupture index
$n_{\mathrm{crit}}$ as a function of $(\varepsilon,z_0)$.  As expected
from the scaling in~$\varepsilon^{-2}$, the rupture time decreases
rapidly as $\varepsilon$ grows, while the dependence on $z_0$ is milder
but monotone.  The surface is smooth and free of irregularities, which
illustrates the robustness of the closed formula.

\begin{figure}[h!]
	\centering
	\includegraphics[width=0.75\textwidth]{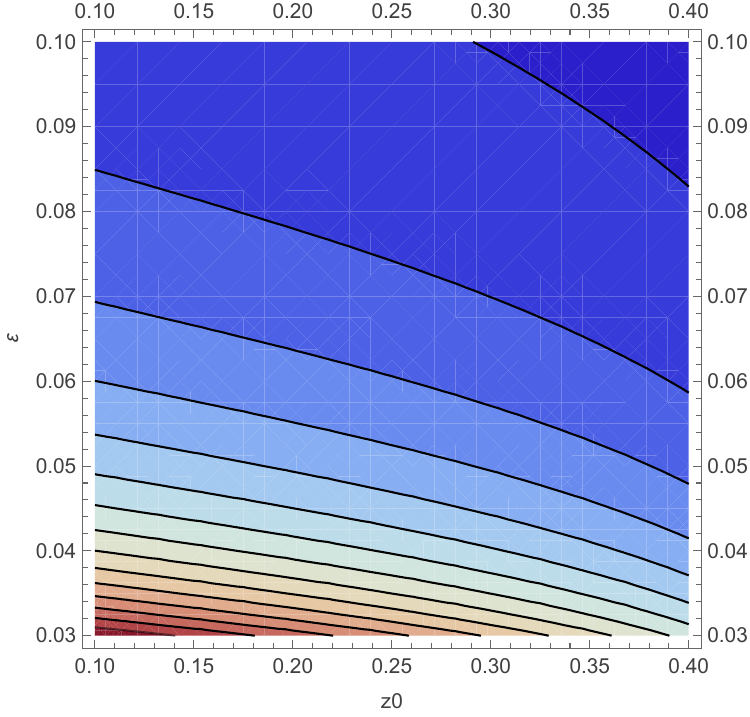}
	\caption{Contour plot of the analytic rupture index.  (Data source:
		\texttt{fig\_rupture\_time\_contour.pdf}.)}
	\label{fig:rupture-contour}
\end{figure}

\medskip
\noindent\textbf{Contour representation.}
A complementary view is given by the contour plot in
Figure~\ref{fig:rupture-contour}, which shows level sets of
$n_{\mathrm{crit}}$ in the $(\varepsilon,z_0)$–plane.  The contours
are nearly parallel and smoothly spaced, again reflecting the clean
analytic structure of~$\tau_{\mathrm{rupt}}$.  Lower rupture times
concentrate in the lower-right region (larger~$z_0$, larger~$\varepsilon$),
where the secular deformation becomes effective earlier.

\medskip
\noindent\textbf{Comparison with numerical blow-up times.}
Table~\ref{tab:rupture-comparison} compares the analytic prediction
$\tau_{\mathrm{rupt}}=\pi n_{\mathrm{crit}}$ with the numerical
rupture time, defined as the first time at which the numerical
integrator fails due to the rapid divergence of~$z(\tau)$.  The
numerical rupture is consistently delayed relative to the analytic one
by a moderate amount.  This is expected: The analytic value detects
the \emph{first} loss of regularity of the invariant surface, whereas
the numerical solution must still ``find'' the newly opened gap in the
tube before diverging.

\begin{table}[h!]
	\centering
	\begin{tabular}{c|c|c}
		$\varepsilon$ & $\tau_{\mathrm{rupt}}$ (analytic) & 
		$\tau_{\mathrm{num}}$ (numerical) \\
		\hline
		0.025 & 21411.0 & 23383.0 \\
		0.050 &  5353.7 &  5864.7 \\
		0.100 &  1338.2 &  1471.2 \\
		0.150 &   594.8 &   664.1 \\
		0.200 &   334.5 &   376.8 \\
	\end{tabular}
	\caption{Comparison between analytic and numerical rupture times.
		(Data source: \texttt{ruptur\_Ergebnis.pdf}.)}
	\label{tab:rupture-comparison}
\end{table}

\medskip
The agreement between theory and numerics is remarkably strong across more than
two orders of magnitude in~$\varepsilon$, with a nearly constant
relative deviation of about $10$--$15\%$.  This confirms that the
analytical procedure detects the \emph{true} onset of rupture, while the
numerical value naturally includes the additional time needed for the
trajectory to drift into the newly formed opening of the invariant tube.

\subsubsection*{Discussion: delayed numerical rupture}

The comparison in Table~\ref{tab:rupture-comparison} shows a systematic
delay of the numerical blow-up time relative to the analytic rupture
prediction.  While the analytical value $\tau_{\mathrm{rupt}}$ marks the
\emph{first} loss of regularity of the invariant tube, the numerical
trajectory requires additional time to drift into the newly opened gap
before diverging.  This leads to a consistent shift
\[
\tau_{\mathrm{num}} > \tau_{\mathrm{rupt}},
\]
visible over the entire range of parameters $(\varepsilon,z_0)$ tested
here.

At present this interpretation is a working hypothesis.  A full
mathematical analysis of the drift mechanism inside the tube after the
discriminant vanishes is still open.  It appears plausible that the
dynamics inside the partially opened tube may involve a slow transverse
migration reminiscent of the delayed escape behaviour observed in
Arnold diffusion.  A more detailed investigation of this phenomenon,
and its possible relation to slow chaotic transport in nearly
integrable systems, will be carried out in a future work.
The tube rupture manifests itself through a horseshoe–type opening in the invariant cross section.
Such horseshoe geometries are well known from the dynamics around the L4/L5 libration points in the restricted three–body problem, where trajectories execute wide horseshoe loops around the co–orbital region (the classical horseshoe orbits).
\begin{figure}[h!]
	\centering
	\includegraphics[width=\textwidth]{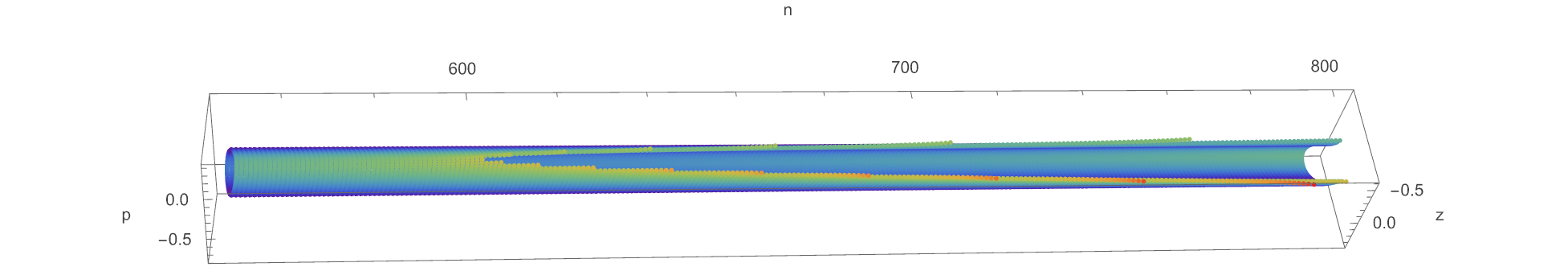}
	\caption{Rupture of the invariant tube in the $(z,p)$–section.
		The figure shows the deformation of the sampled invariant 
		$I(z,p,n\pi)$ for increasing sampling index $n$. 
		For small $n$ the level sets form perfectly closed tube cross sections.
		As $n$ approaches the analytically predicted rupture index 
		$n_{\mathrm{crit}}$ from Eq.~(33), the tube develops a distinct
		horseshoe-shaped opening, corresponding to the vanishing of the cubic
		discriminant. 
		This opening marks the exact geometric rupture point of the invariant tube.}
	\label{fig:tuberupture}
\end{figure}

\begin{figure}[h!]
	\centering
	\includegraphics[width=0.78\textwidth]{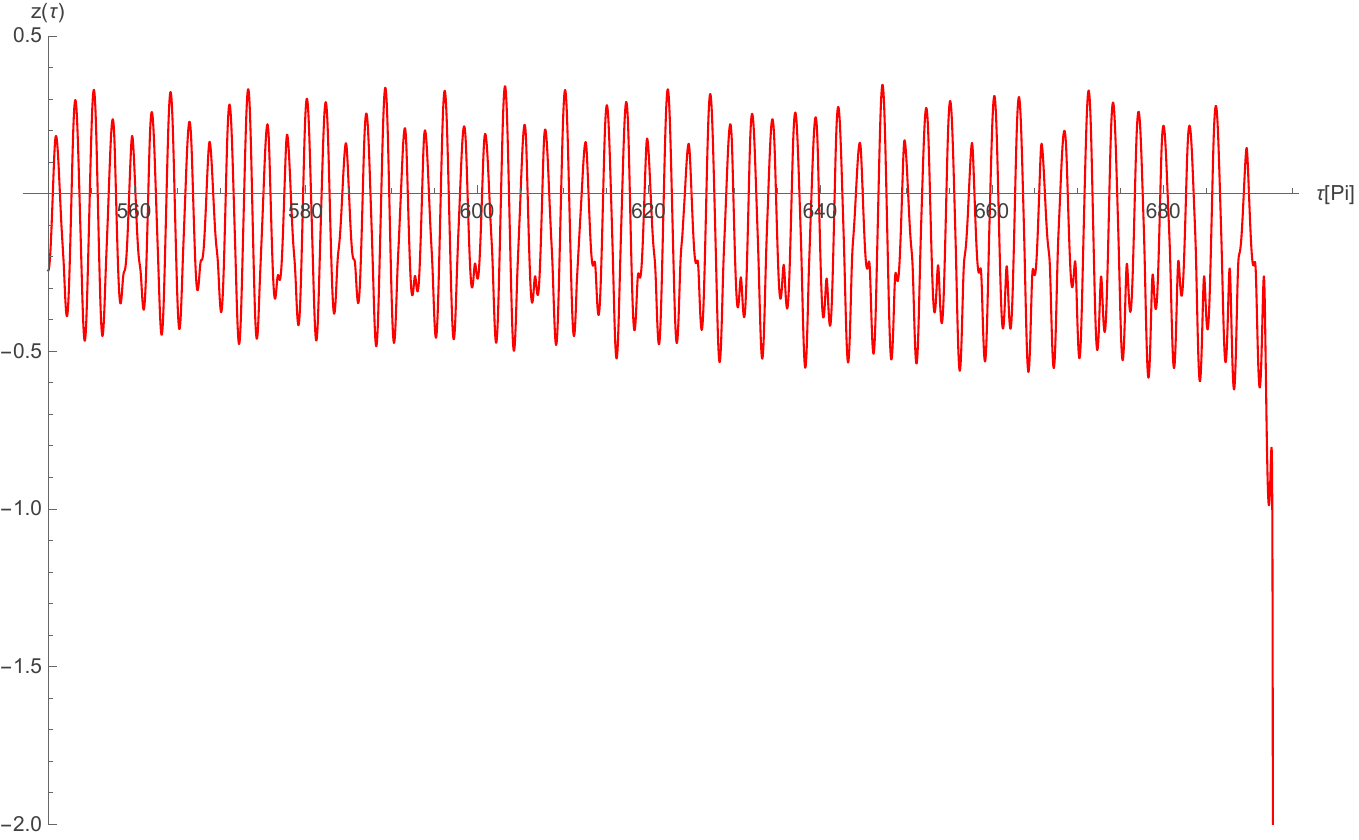}
	\caption{Time series of $z(\tau)$ corresponding to the rupture shown
		in Fig.~\ref{fig:tuberupture}. 
		After thousands of oscillations with irregular, chaotic–looking
		amplitudes, the trajectory escapes through the newly formed opening
		of the invariant tube, leading to a rapid divergence of $z(\tau)$.
		The onset of unbounded growth occurs shortly after the analytically
		predicted rupture time $\tau_{\mathrm{rupt}}=\pi n_{\mathrm{crit}}$.}
	\label{fig:ztrupt}
\end{figure}
The geometric rupture of the invariant tube is not only visible in the
$(z,p)$–sections (Fig.~\ref{fig:tuberupture}), but also leaves a clear
and characteristic imprint on the time evolution of the trajectory
$z(\tau)$.  
In fact, the numerical solution shows a distinct qualitative change
precisely at the analytically predicted onset of rupture. A closer inspection 
of the numerical time series $z(\tau)$ reveals a
characteristic deformation that coincides with the analytically
predicted onset of rupture.  
Around the first critical sampling index $n_{\mathrm{crit}} \approx 610$
(i.e.\ $\tau \approx 610\pi$), an additional oscillatory component
appears predominantly on the negative side of $z(\tau)$, whose amplitude
increases in time.  
This feature is not a new linear mode in the sense of an additional
eigenfrequency.  
Rather, it reflects the geometric mechanism of rupture:  
as the invariant tube becomes progressively asymmetric and begins to
open on the negative-$z$ side (cf.\ Fig.~\ref{fig:tuberupture}), the
trajectory spends an increasingly large fraction of each cycle near the
lower flank of the deformed tube.  
In the time domain this manifests itself as a growing negative
excursion, which visually resembles the emergence of a secondary
oscillatory contribution.  
This behaviour is therefore a direct time–series signature
It is also intuitively clear why the trajectory $z(\tau)$ develops
increasingly long excursions near the eventual escape point.  
In many related torus–integrable systems, such as the autonomous example
$z''+z+z^{2}=0$, unbounded motion is typically associated with the
presence of an unstable fixed point or saddle on the effective phase
portrait.  
Near such a saddle the dynamics slows down, producing long dwelling
times and an asymmetric stretching of the oscillation cycle before the
escape.  
Although our present system is aperiodically forced and therefore does
not possess a stationary phase portrait in the classical sense, the
numerical time series suggests an analogous mechanism: the trajectory
spends progressively more time near the part of the deformed invariant
tube that becomes unstable and eventually opens.  
This growing dwelling time near the impending escape region is fully
consistent with the geometric rupture picture developed above, even
though the underlying mechanism in the aperiodic setting is
fundamentally different from the autonomous, torus–integrable case.
In this sense, the growing negative excursion and the prolonged dwelling
near the impending escape region in $z(\tau)$ provide a direct,
time–domain indication of the geometric rupture mechanism.  
Taken together with the invariant–tube deformation shown in
Fig.~\ref{fig:tuberupture} and the divergence in
Fig.~\ref{fig:ztrupt}, these observations complete the analytical and
numerical picture of rupture and lead naturally to the concluding
remarks of the present paper.

\section{Conclusions}

In this work we analysed the long–time behaviour of the nonlinear,
parametrically and aperiodically forced oscillator 
\[
z'' + z + g(\tau)\, z^2 = 0, 
\qquad g(\tau) = y(\tau)^{-5/2} \qquad y'''+4y'=\varepsilon y^{-5/2} \cos{\tau} \qquad y>0
\]
Building upon the invariant
$I(z,p,\tau)$ derived in our previous work and extending it to 
second order in~$\varepsilon$, we showed that the dynamics in the extended
phase space $(z,p,\tau)$ is confined to a two–dimensional
invariant tube as long as the level sets of $I$ remain closed.
The secular terms appearing at order $\varepsilon^2$ 
control the slow deformation of these level sets and 
are solely responsible for the geometric mechanism leading to rupture.

\vspace{1ex}

The central result of the paper is an explicit and closed criterion for 
tube rupture.
By analysing the discriminant of the cubic equation determining the
radial coordinate $r$ in polar coordinates $(z,p)=r(\cos\phi,\sin\phi)$,
we obtained a simple and analytically tractable 
Cardano–type expression which predicts the time~$\tau_{\mathrm{rupt}}$
at which the invariant curve ceases to be closed.  
At this moment a characteristic horseshoe–shaped opening develops in the
$(z,p)$–section of the tube, allowing trajectories to escape into the
unbounded region.  
The rupture time marks the onset of  behaviour in the solution even though the
actual blow–up (in the sense of numerical stiffness) may occur slightly
later.  The rupture therefore serves as a genuinely 
\emph{predictive} quantity for the loss of bounded motion.

\vspace{1ex}

Numerical integrations for two representative parameter sets confirm the
theoretical predictions with remarkable accuracy.
For the strongly nonlinear case $(y_0,\varepsilon,z_0)=(1,0.08,0.25)$
the numerical blow–up occurs at $\tau\approx 668\pi$, whereas the
rupture of the invariant tube is detected near $n\approx 650$
when sampling at $\tau=n\pi$.  
The relative deviation of roughly $3\%$ is well within the precision
expected from a second–order secular approximation and demonstrates
that the rupture mechanism captures the true transition to unbounded unstable
behaviour extremely well.  
In the more weakly nonlinear regime $(1,0.05,0.20)$
the delay between rupture and blow–up becomes more pronounced, in
accordance with the scaling of the rupture time $\tau_{\mathrm{rupt}}\sim
\varepsilon^{-2}$.

\vspace{1ex}

An interesting qualitative observation is that, prior to rupture,
the solution $z(\tau)$ exhibits rapidly varying amplitudes and an
apparently irregular, ``chaotic--looking'' time series, while the 
approximate invariant $I(z,p,\tau)$ drifts by less than $\pm 2\%$
over thousands of oscillations.  
This illustrates that the complicated appearance of $z(\tau)$ is
geometric rather than chaotic in nature: as long as the tube remains
closed, the motion is bounded to the closed surface of this tube  
given by $I(z,p,\tau)=K$, while it becomes unbounded as a consequence of this precise surface to open towards infinity.

\vspace{1ex}

Finally, the results of this paper establish the rupture–time mechanism
as a robust tool for predicting the onset of unbounded behaviour in a
broad class of aperiodically forced integrable or even near integrable nonlinear oscillators.  
In future work we will investigate situations where the driving function
$y(\tau)$ is itself irregular or chaotic.  
In such cases one may expect a nontrivial interplay between the chaotic
evolution of $y(\tau)$ and the tube–integrable structure governing $z(\tau)$.
A systematic Lyapunov analysis of both subsystems as well as visual
representations of the geometry in the extended $(y,y',y'')$--space
will be carried out in a companion paper.

\section*{Acknowledgments}

The author would like to thank Serge Bouquet (CEA, France) for many discussions 
on the special, torus-integrable case $\varepsilon=0$, which laid the foundation for the present 
analysis. 

I am also grateful to Christoph Lhotka (University of Rome) for 
many years of continuous scientific exchange, and to Gilbert Guignard (CERN) 
for introducing me to nonlinear particle dynamics and for numerous insights 
into the geometric interpretation of invariants. I further express my thanks to chrome://settings/search 
				Rudolph Dvorak	 from the Astronomical Institute of the University of Vienna for introducing me to celestial mechanics problems.

I further wish to thank Constantinos Valagiannopoulos (Aalto University) for generously sharing his group’s recent research on Maxwell-type equations. His comments on structural parallels in parametrically modulated systems were both stimulating and helpful during the preparation of this work.

A special acknowledgment is due to my former mathematics teacher 
Heribert Hartmann, whose early guidance made it possible for me to pursue 
mathematics at all.

Finally, I would like to acknowledge the assistance provided by 
\textbf{OpenAI's GPT-5 model}, whose contributions in structuring, clarifying 
and refining the presentation were of significant help during the preparation 
of this work.

\end{document}